\documentclass[12pt]{article}
\usepackage{amsmath,amssymb,amsfonts}
\usepackage{color}
\numberwithin{equation}{section}
\newtheorem{thm}{Theorem}

\newtheorem{prop}[thm]{Proposition}
\newtheorem{lem}[thm]{Lemma}
\newtheorem{cor}[thm]{Corollary}

{ \pagestyle{plain}

\def\cwedge{\bigcirc\kern-1.07em\wedge\ }

\newcommand{\qed}{\hfill\fbox{}\par\vspace{.2cm}}

\begin{document}

%\title{\bf{{Curvature identities on contact metric manifolds and its applications}} }
%\author{Jeong Hyeong Park$^{*}$, K. Sekigawa$^{\dag}$ and Won Min Shin
%\thanks{Department of Mathematics, Sungkyunkwan University,
% Suwon 440-746, Korea, e-mail: parkj@skku.edu} \thanks{Department of Mathematics,
%    Niigata University,
%    Niigata 950-2181, JAPAN, e-mail: sekigawa@math.sc.niigata-u.ac.jp}
% }
%
%\date{}
%
%\maketitle

%%%%%%%%%%%%%%%%%%%%%%%%%%%%%%%%%%%%%
\begin{center}
{\LARGE {{Curvature identities on contact metric manifolds \\and
their applications}}}
\end{center}
{
\begin{center}
{\large JeongHyeong Park$^{\dag} $, Kouei Sekigawa$^{\ddag}$ and Won
Min Shin$^{\dag}$}
\end{center} }
\begin{center}
$^{\dag}$Sungkyunkwan University,
 Suwon, Korea\\
$^{\ddag}$Niigata University,
    Niigata, Japan\end{center}

%%%%%%%%%%%%%%%%%%%%%%%%%%%%%%%%%%%%%%%%%%%%%%%%%%%%%%%%%%%%%%%%%%%%%%%%%
\begin{abstract}
{{We study curvature identities on contact metric manifolds
           on the geometry of the corresponding almost K\"aehler
           cones,  and we provide applications of the derived curvature
           identities.}}
\end{abstract}
\noindent {\it Mathematics Subsect Classification (2010)} : 53B20, 53C25

%%%%%%%%%%%%%%%%%%%%%%%%%%%%%%%%%%%%%%%%%%%%%%%%%%%%%%%%%%%%%%%%%%%%%%%%%%%%%%%%%
\section{Introduction}\label{sec1}
%%%%%%%%%%%%%%%%%%%%%%%%%%%%%%%%%%%%%%%%%%%%%%%%%%%%%%%%%%%%%%%%%%%%%%%%%%%%%%%%%

A (2n+1)-dimensional smooth manifold $M$ is called an almost contact
manifold if it admits a triplet $(\phi, \xi, \eta)$ of a (1,1)
tensor field $\phi$, a vector field $\xi$ and a 1-form $\eta$
satisfying
\begin{equation}\label{11}
\begin{split}
\phi^2X=-X+\eta(X)\xi,\quad \phi\xi=0, \quad \eta\circ\phi=0, \quad
\eta(\xi)=1,
\end{split}
\end{equation}
for any $X \in\mathfrak{X}(M)$, where $\mathfrak{X}(M)$ denote the
Lie algebra of all smooth vector fields on $M$. Further, an almost
contact metric manifold $M=(M, \phi, \xi, \eta)$ equipped with a
Riemannian metric $g$ satisfying
\begin{equation}\label{12}
\begin{split}
g(\phi{X},\phi{Y})=g(X,Y)-\eta(X)\eta(Y), \quad \eta(X)=g(\xi,X),
\end{split}
\end{equation}
for any $X, Y \in\mathfrak{X}(M)$ is called an almost contact metric
manifold with the almost contact metric structure
$(\phi,\xi,\eta,g)$. On the other hand, a (2n+1)-dimensional smooth
manifold $M$ is called a contact manifold if it admits a global
1-form $\eta$ satisfying $\eta \wedge (d\eta)^n \neq 0$ everywhere
on $M$. Then, the 1-form $\eta$ is called a contact form of $M$. It
is well-known that, for a given contact manifold $M=(M,\eta)$, there
exists an almost contact metric structure $(\phi,\xi,\eta,g)$ ($\xi$
being the dual vector field of $\eta$) satisfying
\begin{equation}\label{13}
\begin{split}
d\eta(X,Y)=g(X,\phi{Y}), \quad g(X,\xi)=\eta(X),
\end{split}
\end{equation}
for any $X, Y \in\mathfrak{X}(M)$, which is called a contact metric
structure on $M$ associated to the contact structure $\eta$. A
contact manifold endowed with an associated contact metric structure
is called a contact metric manifold. Now, let
$M=(M,\phi,\xi,\eta,g$) be a (2n+1)-dimensional almost contact
metric manifold and $\bar{M}=M\times \mathbb R$ be the product
manifold of $M$ and a real line $\mathbb R$ equipped with the
(1,1)-tensor field $\bar{J}$ and a Riemannian manifold $\bar{g}$
defined respectively by
\begin{equation}\label{14}
\begin{split}
&\bar{J}X=\phi{X}-\eta(X)\frac{\partial}{\partial{t}}, \quad
\bar{J}\frac{\partial}{\partial{t}}=\xi,\\
&\bar{g}(X,Y)=e^{-2t}g(X,Y), \quad
\bar{g}(\frac{\partial}{\partial{t}},\frac{\partial}{\partial{t}})=e^{-2t},
\quad \bar{g}(X,\frac{\partial}{\partial{t}})=0,
\end{split}
\end{equation}
for any $X, Y \in\mathfrak{X}(M)$ and $t \in \mathbb R$. Then, we
may easily check that $(\bar{J},\bar{g})$ is an almost Hermitian
structure on $\bar{M}$. An almost contact metric manifold $M=(M,
\phi, \xi, \eta, g)$ is said to be {\it{normal}}  if the
corresponding almost Hermitian manifold $\bar{M}=(\bar{M}, \bar{J},
\bar{g})$ is integrable (i.e. $\bar{M}=(\bar{M}, \bar{J}, \bar{g})$
is a Hermitian manifold). Further, a normal contact metric manifold
is called a Sasakian manifold. Concerning the relationships between
the classes of almost contact metric manifolds and the $cones$ of
the corresponding almost Hermitian manifolds defined by \eqref{14},
it is known that an almost contact metric manifold $M=(M, \phi, \xi,
\eta, g)$ is a Sasakian manifold (resp. a contact metric manifold)
if and only if the corresponding almost Hermitian manifold
$\bar{M}=(\bar{M}, \bar{J}, \bar{g})$ is a $K\ddot{a}hler$ manifold
(resp. an {\it{almost}} $K\ddot{a}hler$ manifold) {\cite{T}}. In
\cite{G}, Gray established curvature identities for an almost
Hermitian manifold belonging to some special classes of almost
Hermitian manifolds, for example, $K\ddot{a}hler$ manifolds,
{\it{almost}} $K\ddot{a}hler$ manifolds, {\it{quasi}}
$K\ddot{a}hler$ manifolds and Hermitian manifolds. In the present
paper, we shall discuss the curvature identities on contact metric
manifolds $M=(M, \phi, \xi, \eta, g)$ derived from the curvature
identities on the corresponding almost Hermitian manifolds
$\bar{M}=(\bar{M}, \bar{J}, \bar{g})$ defined by \eqref{14}. We also
provide some results related to the obtained curvature identities.
%%%%%%%%%%%%%%%%%%%%%%%%%%%%%%%%%%%%%%%%%%%%%%%%%%%%%%%%%%%%%%%%%%%%%%%%%%%%%%%%%%%%%%%
\section{Preliminaries}\label{sec2}
%%%%%%%%%%%%%%%%%%%%%%%%%%%%%%%%%%%%%%%%%%%%%%%%%%%%%%%%%%%%%%%%%%%%%%%%%%%%%%%%%%%%%%%
%{\textbf{2.1. ~3-dimensional contact metric manifolds}}
\indent In this section, we prepare some basic formulas and
fundamental facts we need in the discussions of the present paper.
Let $M=(M,\phi,\xi,\eta,g)$ be a (2n+1)-dimensional almost contact
metric manifold and $\bar{M}=(\bar{M},\bar{J},\bar{g})$ be the
corresponding almost Hermitian manifold defined by \eqref{14}. We
denote by $\nabla$(resp.$\bar{\nabla}$) the Levi-Civita connection
of the Riemannian metric $g$(resp.$\bar{g}$) and by
$R$(resp.$\bar{R}$) the curvature tensor of
$\nabla$(resp.$\bar{\nabla}$) defined respectively by
\begin{equation}\label{21}
\begin{split}
R(X,Y)Z=[\nabla_X, \nabla_Y]Z- \nabla_{[X,Y]}Z,
\end{split}
\end{equation}
for $X, Y, Z \in\mathfrak{X}(M)$, and
\begin{equation}\label{22}
\begin{split}
\bar{R}(\bar{X},\bar{Y})\bar{Z}=[\bar{\nabla}_{\bar{X}},
\bar{\nabla}_{\bar{Y}}]\bar{Z}-
\bar{\nabla}_{[\bar{X},\bar{Y}]}\bar{Z},
\end{split}
\end{equation}
for $\bar{X}, \bar{Y}, \bar{Z} \in\mathfrak{X}(\bar{M})$. We note
that $\mathfrak{X}(\bar{M})$ can be regarded as a Lie subalgebra of
$\mathfrak{X}(\bar{M})$ in the natural way. Further, we get
\begin{equation}\label{23}
\begin{split}
R(X,Y,Z,W)=g(R(X,Y)Z,W)~(resp.~\bar R(\bar X,\bar Y,\bar Z,\bar
W)=\bar g(\bar R(\bar X,\bar Y)\bar Z,\bar W)),
\end{split}
\end{equation}
for ${X}, {Y}, {Z}, W \in\mathfrak{X}({M})$ (resp. for $\bar{X},
\bar{Y}, \bar{Z}, \bar{W} \in\mathfrak{X}(\bar{M})$). We denote also
by $\rho$ and $\tau$ (resp. $\bar{\rho}$ and $\bar{\tau}$) the Ricci
tensor and the scalar curvature of $M=(M,\phi,\xi,\eta,g)$ (resp.
$\bar{M}=(\bar{M},\bar{J},\bar{g})$). Further, we denote by
$\bar{\rho}^*$ and $\bar{\tau}^*$ (resp. $\rho^*$ and $\tau^*$) the
$*$-Ricci tensor and $*$-scalar curvature of
$\bar{M}=(\bar{M},\phi,\xi,\eta,\bar{g})$(resp.
$M=(M,\phi,\xi,\eta,g)$) defined by
\begin{equation}\label{24}
\begin{split}
&\bar{\rho}^*(\bar{X},\bar{Y})=\frac{1}{2}trace \hspace{1mm} of
(\bar{Z}\mapsto
\bar{R}(\bar{X},\bar{J}\bar{Y})\bar{J}\bar{Z}),\\
(resp.~ &\rho^*(X, Y)=\frac{1}{2}trace \hspace{1mm}of (Z \mapsto
R(X,\phi{Y})\phi{Z})),
\end{split}
\end{equation}
and
\begin{equation}\label{25}
\begin{split}
&\bar{\tau}^* = trace \hspace{1mm}of\hspace{1mm}(1,1)~ \text{tensor field}~ \bar{Q}^*~{\text on}~ \bar{M},\\
(resp.~ &\tau^* = trace\hspace{1mm} of\hspace{1mm}(1,1)~\text{tensor field}~ Q^*~{\text on}~ M),\\
\end{split}
\end{equation}
given by
\begin{equation}\label{26}
\begin{split}
&\bar{g}(\bar{Q}^*\bar{X},\bar{Y})=\bar{\rho}^*(\bar{X},\bar{Y}),\\
(resp.~ &g(Q^*X, Y)=\rho^*(X, Y)),\\
\end{split}
\end{equation}
for $\bar{X}, \bar{Y}, \bar{Z} \in\mathfrak{X}(\bar{M})$ $(resp.~X,
Y, Z \in \mathfrak{X}(M))$. Now we define (1,1)-tensor field $h$ by
\begin{equation}\label{27}
\begin{split}
h=\frac{1}{2}\pounds_\xi\phi.
\end{split}
\end{equation}
The tensor field $h$ plays an important role in the geometry of
almost contact metric manifold. The tensor field $h$
satisfies the following equalities
\begin{equation}\label{28}
\begin{split}
h\xi=0, \quad trh=0.
\end{split}
\end{equation}
Further, it is known that the tensor field $h$ is symmetric and
satisfies
\begin{equation}\label{29}
\begin{split}
\phi{h}+h\phi=0,
\end{split}
\end{equation}
for a contact metric manifold $M=(M,\phi,\xi,\eta,g)$. Now, from
\eqref{14}, by direct calculation, we have
\begin{equation}\label{210}
\begin{split}
\bar{\nabla}_XY=\nabla_XY+g(X,Y)\frac{\partial}{\partial{t}}, \quad
\bar{\nabla}_X\frac{\partial}{\partial{t}}=-X, \quad
\bar{\nabla}_\frac{\partial}{\partial{t}}X=-X, \quad
\bar{\nabla}_\frac{\partial}{\partial{t}}\frac{\partial}{\partial{t}}=-\frac{\partial}{\partial{t}},
\end{split}
\end{equation}
for $X, Y \in\mathfrak{X}(M)$. Thus, from \eqref{14} and
\eqref{210}, we have further
\begin{equation}\label{211}
\begin{split}
(\bar{\nabla}_X\bar{J})Y=(\nabla_X\phi)Y-g(X,Y)\xi+\eta(Y)X-(g(\phi{X},Y)+(\nabla_X\eta)(Y))\frac{\partial}{\partial{t}},
\end{split}
\end{equation}
\begin{equation}\label{212}
\begin{split}
(\bar{\nabla}_X\bar{J})\frac{\partial}{\partial{t}}=\nabla_X\xi+\phi{X},
\end{split}
\end{equation}
\begin{equation}\label{213}
\begin{split}
(\bar{\nabla}_\frac{\partial}{\partial{t}}\bar{J})X=0,\quad
(\bar{\nabla}_\frac{\partial}{\partial{t}}\bar{J})\frac{\partial}{\partial{t}}=0,
\end{split}
\end{equation}
for $X, Y \in\mathfrak{X}(M)$. Then, from \eqref{14} and
\eqref{211}, we have
\begin{equation}\label{214}
\begin{split}
\bar{\nabla}_i\bar{J}_{jk}=\nabla_i\phi_{jk}-g_{ij}\eta_k+\eta_jg_{ik}.
\end{split}
\end{equation}
Throughout this paper, we shall adopt the notational convention in
the usual tensor analysis. Let $(x^i)=(x^1,x^2,\ldots,x^{2n+1})$ be
local coordinates on an open subset of $M$ and
$(x^\lambda)=(x^1,x^2,\ldots,x^{2n+1},x^{2n+2}=t)=(x^i,x^{2n+2})$ be
the corresponding local coordinates an open subset $U\times \mathbb
R$ of $\bar{M}=M \times \mathbb R$. Further, we assume that the
Latin indices run over the range $1,2,\ldots , 2n+1$ and the Greek
indices run over the range $1, 2, \ldots , 2n+1, 2n+2= \varDelta$,
and
\begin{equation}\label{215}
\begin{split}
&\phi{\partial_j}={\phi_j}^i\partial_i \quad
(\partial_j=\frac{\partial}{\partial{x^j}}), \quad
\phi_{ij}=g(\phi{\partial_i},\partial_j)=g_{aj}{\phi_i}^a,\\
(resp.~&{\bar{J}\partial_\lambda = \bar{J}_\lambda}^\mu \partial_\mu
~(\partial_\lambda=\frac{\partial}{\partial{x^\lambda}},~
\partial_{\varDelta}=\frac{\partial}{\partial{x^{\varDelta}}}=\frac{\partial}{\partial{t}})\\
&\bar{J}_{\lambda\mu}=\bar{g}(\bar{J}\partial_\lambda,
\partial_\mu)={\bar{g}_{\alpha\mu} \bar{J}_{\lambda}}^\alpha).
\end{split}
\end{equation}
Then, from \eqref{14} and \eqref{215}, we see that \eqref{211},
\eqref{212} and \eqref{213} can be rewritten respectively as
follows:
\begin{equation}\label{216}
\begin{split}
&{\bar{\nabla}_i\bar{J}_j}^k= \nabla_i{\phi_j}^k- g_{ij}\xi^k+\eta_j{\delta_i}^k,\\
&\bar{\nabla}_i{\bar{J}_j}^\varDelta =-\phi_{ij}-\nabla_i\eta_j,\\
&{\bar{\nabla}_i\bar{J}_{\varDelta}}^k=\nabla_i\xi^k+{\phi_i}^k,\\
&{\bar{\nabla}_{\varDelta}\bar{J}_j}^k=0,~{\bar{\nabla}_{\varDelta}\bar{J}_j}^{\varDelta}=0,~{\bar{\nabla}_{\varDelta}\bar{J}_{\varDelta}}^{\varDelta}=0.
\end{split}
\end{equation}
Further, we also set
\begin{equation}\label{217}
\begin{split}
&R(\partial_i,\partial_j)\partial_k={R_{ijk}}^l\partial_l,~R_{ijkl}=g(R(\partial_i,\partial_j)\partial_k,\partial_l)={g_{dl}R_{ijk}}^d,\\
&(resp.~{\bar{R}(\partial_\lambda,\partial_\mu)\partial_\nu=
\bar{R}_{\lambda\mu\nu}}^{\kappa}\partial_k,~\bar{R}_{\lambda\mu\nu
k}=\bar{g}(\bar{R}(\partial_\lambda,\partial_\mu)\partial_\nu,\partial_k)={\bar{g}_{k\sigma}\bar{R}_{\lambda\mu\nu}}^\sigma),
\end{split}
\end{equation}
and so on. Then, from \eqref{14} and \eqref{217}, we have
\begin{equation}\label{218}
\begin{split}
&{{\bar{R}}_{ijk}}~^l={R_{ijk}}^l-{\delta_i}^lg_{jk}+{\delta_j}^lg_{ik},\\
&{\bar{R}_{ijk}}~^{\varDelta}=0, \quad {\bar{R}_{i\varDelta
k}}~^l=0, \quad{\bar{R}_{i\varDelta\varDelta}}~^l=0, \quad
{\bar{R}_{i\varDelta\varDelta}}~^{\varDelta}=0.
\end{split}
\end{equation}
From \eqref{218}, we also have
\begin{equation}\label{219}
\begin{split}
\bar{\rho}_{jk}=\rho_{jk}-2ng_{jk},\quad \bar{\rho}_{j \varDelta}=0,
\quad \bar{\rho}_{\varDelta\varDelta}=0,
\end{split}
\end{equation}
and hence
\begin{equation}\label{220}
\begin{split}
\bar{\tau}=e^{2t}(\tau-4n^2-2n).
\end{split}
\end{equation}
Similarly, from \eqref{14}, \eqref{24}, and \eqref{218}, we have
\begin{equation}\label{221}
\begin{split}
\bar{\rho}^*_{ij}=\rho^*_{ij}-g_{ij}+\eta_i\eta_j,\quad
\bar{\rho}^*_{i\varDelta}=\frac{1}{2}\xi^a{R_{iac}}^b{\phi_b}^c,
\quad \bar{\rho}^*_{\varDelta j}=0, \quad
\bar{\rho}^*_{\varDelta\varDelta}=0,
\end{split}
\end{equation}
and hence
\begin{equation}\label{222}
\begin{split}
\bar{\tau}^*=e^{2t}(\tau^*-2n).
\end{split}
\end{equation}
From \eqref{14} and \eqref{219}, we get
\begin{equation}\label{223}
\begin{split}
&\bar{J}_{j}^\alpha{\bar{J}_k}^\beta\bar{\rho}_{\alpha\beta}=\rho_{ab}{\phi_j}^a{\phi_k}^b-2{n}g_{jk}+2{n}\eta_j\eta_k,\\
&{\bar{J}_j}^\alpha{\bar{J}_{\varDelta}}^\beta\bar{\rho}_{\alpha\beta}=\rho_{ab}{\phi_j}^a\xi^b,\\
&{\bar{J}_{\varDelta}}^\alpha{\bar{J}_{\varDelta}}^\beta\bar{\rho}_{\alpha\beta}=\rho_{ab}\xi^a\xi^b-2{
n}.
\end{split}
\end{equation}
Similarly, from \eqref{14} and \eqref{221}, we get also
\begin{equation}\label{224}
\begin{split}
&{{\bar{J}_j}^\alpha{\bar{J}_k}^\beta\bar{\rho}_{\alpha\beta}^*=\rho_{ab}^*{\phi_j}^a{\phi_k}^b-g_{jk}+\eta_j\eta_k+\rho_{lj}^*\xi^l\eta_k,}\\
&{{\bar{J}_{\varDelta}}^\alpha{\bar{J}_k}^\beta\bar{\rho}_{\alpha\beta}^*=\frac{1}{2}\xi^a{R_{kac}}^b{\phi_b}^c,}\\
&{{\bar{J}_j}^\alpha{\bar{J}_{\varDelta}}^\beta\bar{\rho}_{\alpha\beta}^*={\bar{J}_{\varDelta}}^\alpha{\bar{J}_{\varDelta}}^\beta\bar{\rho}_{\alpha\beta}^*=0.}
\end{split}
\end{equation}
\noindent Since
{{${\bar{J}_{\lambda}}^\alpha{\bar{J}_{\mu}}^\beta\bar{\rho}^*_{\alpha\beta}={\bar{\rho}^*}_{{\mu}{\lambda}}$}}
holds, from \eqref{221} and \eqref{224}, we have
\begin{equation}\label{225}
\begin{split}
\rho^*_{ab}{\phi_j}^a{\phi_k}^b+\rho^*_{bj}\eta_k\xi^b=\rho^*_{kj}.
\end{split}
\end{equation}
Further, from \eqref{14} and \eqref{218}, by direct calculation, we
have
\begin{equation}\label{226}
\begin{split}
&(1)~e^{2t}{\bar{J}_i}^\alpha{\bar{J}_j}^\beta{\bar{J}_k}^\gamma{\bar{J}_l}^\sigma\bar{R}_{\alpha\beta\gamma\sigma}
={\phi_i}^a{\phi_j}^b{\phi_k}^c{\phi_l}^dR_{abcd}-g_{il}g_{jk}+g_{ik}g_{jl}+g_{il}\eta_j\eta_k\\
&\qquad\qquad\qquad\qquad\qquad\qquad+g_{jk}\eta_i\eta_l-{{g_{jl}\eta_i\eta_k}} -g_{ik}\eta_j\eta_l,\\
&(2)~e^{2t}{\bar{J}_i}^\alpha{\bar{J}_j}^\beta\bar{R}_{\alpha\beta
kl}={\phi_i}^a{\phi_j}^bR_{abkl}-\phi_{il}{{\phi}}_{jk}+{{\phi}}_{jl}\phi_{ik},\\
&(3)~e^{2t}{\bar{J}_i}^\alpha{\bar{J}_k}^\gamma{\bar{R}}_{\alpha{{j}}\gamma
l}={\phi_i}^a{\phi_k}^cR_{ajcl}+\phi_{il}\phi_{jk}+g_{jl}g_{ik}-g_{jl}\eta_i\eta_k,\\
&(4)~e^{2t}{\bar{J}_{\varDelta}}^\alpha{\bar{J}_j}^\beta\bar{R}_{\alpha\beta
kl}=\xi^a{\phi_j}^bR_{abkl}-\eta_l\phi_{jk}+\eta_k\phi_{jl},\\
&(5)~e^{2t}{\bar{J}_{\varDelta}}^\alpha{\bar{J}_k}^\gamma\bar{R}_{\alpha
j\gamma l}={\xi}^a{\phi_k}^cR_{ajcl}+\eta_l\phi_{jk},\\
&(6)~e^{2t}{\bar{J}_{\varDelta}}^\alpha{\bar{J}_j}^\beta{\bar{J}_k}^\gamma{\bar{J}_l}^\sigma\bar{R}_{\alpha\beta
\gamma\sigma}=\xi^a{\phi_j}^b{\phi_k}^c{\phi_l}^dR_{abcd},\\
&(7)~e^{2t}{\bar{J}_{\varDelta}}^\alpha{\bar{J}_{\varDelta}}^\gamma\bar{R}_{\alpha
j \gamma l}=\xi^a\xi^cR_{ajcl}+g_{jl}-\eta_j\eta_l,\\
&(8)~e^{2t}{\bar{J}_{\varDelta}}^\alpha{\bar{J}_j}^\beta{\bar{J}_{\varDelta}}^\gamma{\bar{J}_l}^\sigma\bar{R}_{\alpha\beta
\gamma\sigma}=\xi^a{\phi_j}^b\xi^c{\phi_l}^dR_{abcd}+g_{jl}-\eta_j\eta_l,
\end{split}
\end{equation}
%%%%%%%%%%%%%%%%%%%%%%%%%%%%%%%%%%%%%%%%%%%%%%%%%%%%%%%%%%%%%%%%%%%%%%%%%%%%%%%%%%%%%%%
\section{Curvature identities on contact metric manifolds}\label{sec3}
%%%%%%%%%%%%%%%%%%%%%%%%%%%%%%%%%%%%%%%%%%%%%%%%%%%%%%%%%%%%%%%%%%%%%%%%%%%%%%%%%%%%%%%
Let $M=(M, \phi, \xi, \eta, g)$ be a (2n+1)-dimensional contact
metric manifold and $\bar{M}=M \times \mathbb R$ be the product
manifold $M$ and a real line $\mathbb R$ endowed with the almost
Hermitian {{structure}} $(\bar{J}, \bar{g})$ defined by \eqref{14}.
Now, we suppose that the curvature tensor $\bar{R}$ of
        $\bar{M}$ satisfies the following identity:  %\cite{G}:
\begin{equation}\label{31}
\begin{split}
&\bar{R}_{\lambda\mu\nu\kappa}
+{\bar{J}_\lambda}^\alpha{\bar{J}_\mu}^\beta{\bar{J}_\nu}^\gamma{\bar{J}_\kappa}^\sigma
\bar{R}_{\alpha\beta\gamma\sigma}
+{\bar{J}_\lambda}^\alpha{\bar{J}_\nu}^\gamma\bar{R}_{\alpha\mu\gamma\kappa}
+{\bar{J}_\mu}^\beta{\bar{J}_\kappa}^\sigma \bar{R}_{\lambda\beta\nu\sigma} \\
&-{\bar{J}_\lambda}^\alpha{\bar{J}_\mu}^\beta
{{\bar{R}_{\alpha\beta\nu\kappa}}}
-{\bar{J}_\nu}^\gamma{\bar{J}_\kappa}^\sigma
\bar{R}_{\lambda\mu\gamma\sigma}
+{\bar{J}_\lambda}^\alpha{\bar{J}_\kappa}^\sigma
\bar{R}_{\alpha\mu\nu\sigma}
+{\bar{J}_\mu}^\beta{\bar{J}_\nu}^\gamma \bar{R}_{\lambda\beta\gamma\kappa} \\
&=2\bar{g}^{\alpha\beta}(\bar{\nabla}_\alpha
\bar{J}_{\lambda\mu})\bar{\nabla}_\beta\bar{J}_{\nu\kappa}.
\end{split}
\end{equation}
Now, setting $\lambda=i, ~\mu=j, ~\nu=k,~\kappa=l$ in \eqref{31} and
taking account of \eqref{14}, \eqref{215} $\thicksim$
\eqref{218} and \eqref{226}, we have
\begin{equation}\label{32}
\begin{split}
&R_{ijkl}+{\phi_i}^a{\phi_j}^b{\phi_k}^c{\phi_l}^dR_{abcd}+{\phi_i}^a{\phi_k}^cR_{ajcl}+{\phi_j}^b{\phi_l}^dR_{ibkd}-{\phi_i}^a{\phi_j}^bR_{abkl}\\
&-{\phi_k}^c{\phi_l}^dR_{ijcd}+{\phi_i}^a{\phi_l}^dR_{ajkd}
+{\phi_j}^b{\phi_k}^cR_{ibcl}+4\phi_{il}\phi_{jk}-4\phi_{jl}\phi_{ik}\\
&-4g_{il}g_{jk}+4g_{jl}g_{ik}+4g_{il}\eta_j\eta_k+4g_{jk}\eta_i\eta_l-4g_{jl}\eta_i\eta_k-4g_{ik}\eta_j\eta_l\\
&=2(\nabla^a\phi_{ij})\nabla_a\phi_{kl}+2\eta_k\nabla_l\phi_{ij}-2\eta_l\nabla_k\phi_{ij}+2\eta_i\nabla_j\phi_{kl}-2\eta_j\nabla_i\phi_{kl},
\end{split}
\end{equation}
\noindent Similarly, setting $\lambda=\varDelta, ~\mu=j, ~\nu=k,
~\kappa=l$ in \eqref{31}, we have also
%(This is nothing but \eqref{3300})\\
\begin{equation}\label{33}
\begin{split}
&\xi^a{\phi_j}^b{\phi_k}^c{\phi_l}^dR_{abcd}+\xi^a{\phi_k}^cR_{ajcl}-\xi^a{\phi_j}^bR_{abkl}+\xi^a{\phi_l}^dR_{ajkd}\\
&=2(\nabla^a\eta_j)\nabla_a\phi_{kl}+2\eta_k\nabla_l\eta_j-2\eta_l\nabla_k\eta_j+2\phi_{aj}\nabla^a\phi_{kl}.
\end{split}
\end{equation}
\noindent Furthermore, setting $\lambda=i, ~\mu=\varDelta, ~\nu=k,
~\kappa=\varDelta$ in \eqref{31}, we have
\begin{equation}\label{34}
\begin{split}
&{\phi_i}^a\xi^b{\phi_k}^c\xi^dR_{abcd}+\xi^b\xi^dR_{ibkd}
=2((\nabla_a\eta_{{i}})\nabla^a\eta_k+\phi_{ak}\nabla^a\eta_i+\phi_{ai}\nabla^a\eta_k),
\end{split}
\end{equation}
where $\nabla^a=g^{ab}\nabla_b$. Thus, from \eqref{32} $\thicksim$
\eqref{34}, taking account of the $\mathbb Z_2$-symmetricity of the
curvature tensors $R$ and $\bar{R}$ of type
       (0,4) of $M$ and $\bar{M}$ respectively, we have the following.
   \begin{thm}\label{Th1}
The curvature tensor $R$ of an almost contact metric manifold $M=(M,
\phi, \xi, \eta, g)$ verifies the identities \eqref{32} $\thicksim$
\eqref{34} if and only if the curvature tensor $\bar R$ of the
corresponding almost Hermitian manifold $\bar M=(\bar M, \bar J,
\bar g)$ defined by \eqref{14} verifies \eqref{31}.
\end{thm}
Particularly, if  $M$ is a contact metric manifold, then the
corresponding
       almost Hermitian manifold $\bar{M}$ is an {{almost K\"{a}hler}} manifold,
       and hence, the curvature tensor ${\bar{R}}$ of $\bar{M}$ verifies the
       identity \eqref{31}. Thus, from Theorem \ref{Th1}, we have immediately the
       following.
\begin{cor}\label{Co1}
The curvature tensor $R$ of a contact metric manifold
       $M =(M,\phi,\xi,\eta,g)$ verifies the identities \eqref{32} $\thicksim$ \eqref{34}.
\end{cor}
       {{Now, in the remaing of this section}}, we assume that $M = (M,\phi,\xi,\eta,g)$ is a
       (2n+1)-dimensional contact metric manifold and $\bar{M}
       = (\bar{M} \bar{J},{g})$ is the corresponding {{almost
       K\"{a}hler}}
       manifold defined by \eqref{14}.
\noindent We here recall several fundamental formulas on a contact
metric manifold. On the contact metric manifold $M=(M, \phi, \xi,
\eta, g)$ under discussion, we have the following in addition to
\eqref{28} and \eqref{29} \cite{B1}:
\begin{equation}\label{35}
\begin{split}
&(1)~\nabla_i\xi^i=0,\\
&(2)~\xi^a\nabla_a{\phi_j}^i=0,\\
&(3)~\nabla_i\xi^j=-{\phi_i}^j-{\phi_a}^j{h_i}^a,\\
&(4)~\xi^a\nabla_a{h_j}^i={\phi_j}^i-{h_a}^i{h_b}^a{\phi_j}^b-{\phi_a}^i{R_{jbc}}^a\xi^b\xi^c,\\
&{{(5)~\nabla_i{\phi_j}^i=-2n\eta_j,}}\\
&(6)~\rho_{ij}\xi^i\xi^j=2n-trh^2.
\end{split}
\end{equation}
Since
$\bar{\nabla}_i\bar{J}_{jk}+\bar{\nabla}_j\bar{J}_{ki}+\bar{\nabla}_k\bar{J}_{ij}=0$
holds on $M$, from \eqref{214}, we have
\begin{equation}\label{36}
\begin{split}
\nabla_i\phi_{jk}+\nabla_j\phi_{ki}+\nabla_k\phi_{ij}=0.
\end{split}
\end{equation}
From (3) of \eqref{35}, we have also
\begin{equation}\label{37}
\begin{split}
\| \nabla\eta \|^2=2n+trh^2.
\end{split}
\end{equation}
\noindent Transvecting \eqref{32} with $g^{il}$ and taking
accounting of \eqref{215}, (5) of \eqref{35}, we have
\begin{equation}\label{38}
\begin{split}
&2\rho_{jk}+2{\phi_j}^b{\phi_k}^c\rho_{bc}-\xi^a\xi^d{\phi_j}^b{\phi_k}^cR_{abcd}-2\rho^*_{jk}-2\rho^*_{kj}-\xi^a\xi^bR_{ajkb}\\
&-8(n-1)g_{jk}+8(n-1)\eta_j\eta_k=-2(\nabla^a\phi_{jl})\nabla_a{\phi_k}^l+8n\eta_j\eta_k+4\eta_l\nabla_k{\phi_j}^l.
\end{split}
\end{equation}
Further, transvecting \eqref{38} with $g^{jk}$ and taking account of
\eqref{215}, (5), (6) of \eqref{35}, we have
\begin{equation*}
\begin{split}
4\tau-4\tau^*-4(2n-trh^2)-16n^2=-2\| \nabla\phi \|^2,
\end{split}
\end{equation*}
and hence
\begin{equation}\label{39}
\begin{split}
\tau^*-\tau+4n^2=trh^2+\frac{1}{2}(\| \nabla\phi \|^2-4n).
\end{split}
\end{equation} (\cite{O}, (3.8) in Lemma 3.3).
\noindent {Transvecting \eqref{33} with $g^{jl}$ and taking account
of \eqref{215}, (1) of \eqref{35} and \eqref{36}, we have
\begin{equation}\label{310}
\begin{split}
-2\xi^a{\phi_k}^c\rho_{ac}-\xi^a{\phi_j}^b{R_{akb}}^j=2(\nabla_a\eta_j)\nabla^a{\phi_k}^j.
\end{split}
\end{equation}
Furthermore, transvecting \eqref{310} with ${\phi_i}^k$ and taking
account of \eqref{24}, \eqref{28}, \eqref{29}, \eqref{35} and
\eqref{36}, we have
\begin{equation}\label{311}
\begin{split}
\rho_{ai}\xi^a-\rho_{ai}^*\xi^a-2n\eta_i={h_j}^b\nabla_b{\phi_i}^j.
\end{split}
\end{equation}
\noindent On the other hand, transvecting \eqref{31} with
{{$\bar{g}^{\lambda \kappa}$}}, we have
\begin{equation}\label{312}
\begin{split}
\bar\rho_{\mu\nu}^*+\bar\rho_{\nu\mu}^*-\bar\rho_{\mu\nu}-{\bar{J}_\mu}^\alpha{\bar{J}_\nu}^\beta\bar\rho_{\alpha\beta}
={{{\bar{g}^{\alpha\beta}\bar{g}_{\lambda
\sigma}(\bar\nabla_\alpha{\bar{J}_\mu}^\lambda)\bar\nabla_\beta{\bar{J}_\nu}^
\sigma}}}.
\end{split}
\end{equation}
Furthermore, transvecting  \eqref{312} with $\bar{g}^{\mu\nu}$, we
have
\begin{equation}\label{313}
\begin{split}
2(\bar\tau^*-\bar\tau)=\| \bar\nabla\bar{J} \|^2.
\end{split}
\end{equation}
Here, setting $\mu=j,~\nu=k$ in \eqref{312} and taking account of
\eqref{216}, \eqref{219}, \eqref{221}, \eqref{223}, {{(2),}} (5) of
\eqref{35} and \eqref{36}, we have
\begin{equation}\label{314}
\begin{split}
\rho_{jk}^*+\rho_{kj}^*-\rho_{jk}-\rho_{ab}{\phi_j}^a{\phi_k}^b=&(\nabla^a\phi_{jl})\nabla_a{\phi_k}^l-2\eta_l\nabla_j{\phi_k}^l+\phi_{aj}\nabla^a\eta_k+\phi_{ak}\nabla^a\eta_j\\
&+(\nabla^a\eta_j)\nabla_a\eta_k-4(n-1)g_{jk}-4\eta_j\eta_k.
\end{split}
\end{equation}
\noindent Similarly, setting $\mu=i,~\nu=\varDelta$ in \eqref{312},
we have
\begin{equation}\label{315}
\begin{split}
\frac{1}{2}\xi^a{R_{iac}}^b{\phi_b}^c-\rho_{bc}{\phi_i}^b\xi^c=
-(\nabla^a\phi_{ci})\nabla_a\xi^c.
\end{split}
\end{equation}
We may easily check that \eqref{314} can be {{also}} derived from
\eqref{34} and \eqref{38} and further, {{\eqref{315}}} is nothing
but \eqref{310}. So, we can not deduce any new equalities from
{{\eqref{310}}} itself. From \eqref{220} and \eqref{222}, we have
\begin{equation}\label{316}
\begin{split}
{\bar\tau}^*-\bar\tau =e^{2t}(\tau^*-\tau+4n^2).
\end{split}
\end{equation}
On one hand, from \eqref{215}, \eqref{216} and (5) of
\eqref{35}, we have
\begin{equation}\label{317}
\begin{split}
g^{ab}(\bar\nabla_a{\bar{J}_d}^c)\bar\nabla_b{\bar{J}_c}^d&=-(\nabla^a\phi^{cd}-\xi^dg^{ac}+\xi^cg^{ad})(\nabla_a\phi_{cd}-\eta_dg_{ac}+\eta_cg_{ad})\\
&=-(\| \nabla\phi \|^2-4n).
\end{split}
\end{equation}
Thus, taking account of \eqref{317}, we have further
\begin{equation}\label{318}
\begin{split}
\| \bar\nabla\bar{J} \|^2=2(\tau^*-\tau+4n^2)=\| \nabla\phi
\|^2-4n+2trh^2,
\end{split}
\end{equation}
and hence
\begin{equation}\label{319}
\begin{split}
\tau^*-\tau+4n^2=trh^2+\frac{1}{2}(\| \nabla\phi \|^2-4n).
\end{split}
\end{equation}
Thus, we obtain the same equality {{as in (\cite{B1}, Proposition
7.7)}}. It is well-known that a Sasakian manifold is characterized
as an almost contact metric manifold $M=(M, \phi, \xi, \eta, g)$
satisfying the condition} (\cite{B1}, Theorem 6.3):
\begin{equation}\label{320}
\begin{split}
(\nabla_X\phi)Y=g(X, Y)\xi-\eta(X)Y~(namely,
\nabla_i{\phi_j}^k=g_{ij}\xi^k-\eta_i{\delta_j}^k),
\end{split}
\end{equation}
for any X, Y $\in\mathfrak{X}(M)$. Thus summing up the above
arguments, we may reprove the result due to Olszak (\cite{O}, Lemma
3.3).
%%%%%%%%%%%%%%%%%%%%%%%%%%%%%%%%%%%%%%%%%%%%%%%%%%%%%%%%%%%%%%%%%%%%%%%%%%%%%%%%%%%%%%%
\section{Applications} \label{sec4}
%%%%%%%%%%%%%%%%%%%%%%%%%%%%%%%%%%%%%%%%%%%%%%%%%%%%%%%%%%%%%%%%%%%%%%%%%%%%%%%%%%%%%%%
In this section, we shall provide some applications of the
discussions in the previous sections. Gray \cite{G} has defined the
following identities $(\bar{G_i})$ ($i$=1, 2, 3) for the curvature
tensor $\bar{R}$ of an almost Hermitian manifold $\bar M =
(\bar{M},\bar{J},\bar{g})$:
\begin{equation*}
\begin{split}
&(\bar G_1)~ \bar{R}(\bar{X}, \bar{Y}, \bar{Z},
\bar{W})=\bar{R}(\bar{X}, \bar{Y}, \bar{J}\bar{Z},
\bar{J}\bar{W}),\\
&(\bar G_2)~ \bar{R}(\bar{X}, \bar{Y}, \bar{Z},
\bar{W})=\bar{R}(\bar{J}\bar{X}, \bar{J}\bar{Y}, \bar{Z},
\bar{W})+\bar{R}(\bar{J}\bar{X}, \bar{Y}, \bar{J}\bar{Z}, \bar{W})+\bar{R}(\bar{J}\bar{X}, \bar{Y}, \bar{Z}, \bar{J}\bar{W}),\\
&(\bar G_3)~ \bar{R}(\bar{X}, \bar{Y}, \bar{Z},
\bar{W})=\bar{R}(\bar{J}\bar{X}, \bar{J}\bar{Y}, \bar{J}\bar{Z},
\bar{J}\bar{W}),\\
\end{split}
\end{equation*}
\noindent for any $\bar{X}, \bar{Y}, \bar{Z}, \bar{W}
\in\mathfrak{\bar X}(M)$. Almost Hermitian manifold $\bar M =
(\bar{M},\bar{J},\bar{g})$ is said to satisfy the $(\bar{G_i})$
identity if it curvature tensor $\bar{R}$ verifies the curvature
identity $(\bar{G_i})$ ($i$=1, 2, 3). Then we may easily check that
the implication relations for the $(\bar{G_i})$ identities are as
follows \cite{G}:
\begin{equation} \label{41}
~(\bar G_1) \Rightarrow(\bar G_2) \Rightarrow (\bar G_3).
\end{equation}
Now, let $M=(M, \phi, \xi, \eta, g)$ be a (2n+1)-dimensional almost
contact metric manifold and $\bar{M}=M\times \mathbb R$ be the
product manifold of $M$ and a real line $\mathbb R$ endowed with the
almost Hermitian structure $(\bar{J}, \bar{g})$ defined by
\eqref{14}.\\
Now, taking account of \eqref{217}, \eqref{218} and \eqref{226}, we
have the following.
\begin{lem}\label{Le1}
$\bar M = (\bar{M},\bar{J},\bar{g})$ satisfies the
$(\bar{G_1})$ identity if and only if $M=(M, \phi, \xi, \eta, g)$
satisfies
\begin{equation*}
\begin{split}
&(G_1)~ R_{ijkl}-{\phi_k}^c{\phi_l}^dR_{ijcd}=g_{il}g_{jk}-g_{jl}g_{ik}-\phi_{il}\phi_{jk}+\phi_{jl}\phi_{ik},\\
&(G_1-1)~\xi^c{\phi_l}^dR_{ijcd}=\eta_i\phi_{jl}-\eta_j\phi_{il}.
\end{split}
\end{equation*}
\end{lem}
\begin{lem}\label{Le2}
$\bar M = (\bar{M},\bar{J},\bar{g})$ satisfies the
$(\bar{G_2})$ identity if and only if $M=(M, \phi, \xi, \eta, g)$
satisfies
\begin{equation*}
\begin{split}
&(G_2)~ R_{ijkl}-{\phi_i}^a{\phi_j}^bR_{abkl}-{\phi_i}^a{\phi_k}^cR_{ajcl}-{\phi_i}^a{\phi_l}^dR_{ajkd}=g_{jk}\eta_i\eta_l-g_{jl}\eta_i\eta_k,\\
&(G_2-1)~\xi^a{\phi_j}^bR_{abkl}+\xi^a{\phi_k}^cR_{ajcl}+\xi^a{\phi_l}^dR_{ajkd}=0,\\
&(G_2-2)~{\phi_i}^a\xi^bR_{abkl}=\eta_k\phi_{il}-\eta_l\phi_{ik},\\
&(G_2-3)~{\phi_i}^a\xi^cR_{ajcl}=\eta_j\phi_{il},\\
&(G_2-4)~\xi^a\xi^cR_{ajcl}=-g_{jl}+\eta_j\eta_l.
\end{split}
\end{equation*}
\end{lem}
\begin{lem}\label{Le3}
$\bar M = (\bar{M},\bar{J},\bar{g})$ satisfies the
$(\bar{G_3})$ identity if and only if $M=(M, \phi, \xi, \eta, g)$
satisfies
\begin{equation*}
\begin{split}
&(G_3)~ R_{ijkl}-{\phi_i}^a{\phi_j}^b{\phi_k}^c{\phi_l}^dR_{abcd}=g_{il}\eta_j\eta_k+g_{jk}\eta_i\eta_l-g_{ik}\eta_j\eta_l-g_{jl}\eta_i\eta_k,\\
&(G_3-1)~\xi^a{\phi_j}^b{\phi_k}^c{\phi_l}^dR_{abcd}=0,\\
&(G_3-2)~\xi^a\xi^c{\phi_j}^b{\phi_l}^dR_{abcd}=-g_{jl}+\eta_j\eta_l.
\end{split}
\end{equation*}
\end{lem}
  Now, we can easily check that the identity
         $(G_1-1)$  is derived from $(G_1)$ in Lemma \ref{Le1}
         and similarly that the identities
         $(G_2-1) \thicksim (G_2-4)$ are derived from $(G_2)$ in Lemma \ref{Le2},
         and further that the identities $(G_3-1)$
         and $(G_3-2)$ are derived from $(G_3)$ in Lemma \ref{Le3}. Therefore, we have the following.
\begin{prop}\label{Pr1}
$\bar{M} =(\bar{M},\bar{J},\bar{g})$
         satisfies the $(\bar{G}_i)$ identity if and only if
         $M =(M,\phi,\xi,\eta,g)$ satisfies the $(G_i)$
         identity
         $(i = 1,2,3)$.
\end{prop}
\noindent Based on Proposition \ref{Pr1}, Mocanu and Munteanu
\cite{MM}
          defined that an almost contact metric manifold
          $M=(M,\phi,\xi,\eta,g)$ is said to satisfy the $(G_i)$ identity
          if the curvature tensor $R$ of $M$ verifies the identity
          $(G_i)$ ($i$=1, 2, 3).
\noindent Then, from the corresponding implication relations
\eqref{41}, we have the following implication relations for the
$(G_i)$ identities ($i$=1, 2, 3):
\begin{equation}\label{42}
~(G_1) \Rightarrow(G_2) \Rightarrow (G_3).
\end{equation}
The following is well-known (\cite{B1}, Proposition 7.6).
\begin{prop}\label{Pr2}
A contact metric manifold $M=(M, \phi, \xi, \eta, g)$ is Sasakian if
and only if its curvature tensor $R$ satisfies $$R(X,Y)\xi = \eta(Y)
X - \eta(X)Y.$$
\end{prop}
The following is an improvement of the one by Mocanu and
Munreanu \cite{MM}.
\begin{prop}\label{Pr3}
{If {{a}} contact metric manifold $M=(M, \phi, \xi, \eta, g)$
satisfies the $(G_3)$ identity, then $M$ is Sasakian.}
\end{prop}
{\bf{Proof.}} We assume that $M=(M, \phi, \xi, \eta, g)$ is a
contact metric manifold satisfying $(G_3)$ identity. Then,
transvecting the identity $G_3$ with $\xi^k$, we have
\begin{equation*}
\begin{split}
\xi^kR_{ijkl}=g_{il}\eta_j-g_{jl}\eta_i.
\end{split}
\end{equation*}
\noindent Therefore, from Proposition \ref{Pr2}, it follows that $M$ is Sasakian.\\
From Proposition \ref{Pr3}, taking account of \eqref{42}, we have
immediately the following:
\begin{cor}\label{Co2}
{If a contact metric manifold $M=(M, \phi, \xi, \eta, g)$ satisfies
any {{one}} of the $(G_i)$ identities $(i= 1,2,3)$, then M is
Sasakian.}
\end{cor}
In the remaining of {{this}} section, we assume that $M=(M, \phi,
\xi, \eta, g)$ is a 3-dimensional contact metric manifold. Then the
corresponding {\it almost} $K\ddot{a}hler$ manifold $\bar M=(\bar
M,\bar J, \bar g)$ defined by \eqref{14} is a 4-dimensional {\it
almost} $K\ddot{a}hler$ manifold. Since $dim \bar M$=4, it follows
that the following identity.
\begin{equation}\label{43}
\bar{\rho}^*_{\lambda\mu}+\bar{\rho}^*_{\mu\lambda}-(\bar{\rho}_{\lambda\mu}+{{\bar{\rho}_{\alpha\beta}{\bar{J}_\lambda}^{\alpha}{\bar{J}_\mu}^{\beta})}}
=\frac{1}{2}(\bar{\tau}^*-\bar{\tau})\bar{g}_{\lambda\mu},
\end{equation}
holds on $\bar M$ \cite{GH}. Setting $\lambda=i$ and $\mu=j$ in
\eqref{43}, and taking account of \eqref{219} $\thicksim$
\eqref{222}, we have
\begin{equation}\label{43a}
\begin{split}
\rho^*_{ij}+\rho^*_{ji}-\rho_{ij}-\rho_{ab}{\phi_{i}}^{a}{\phi_{j}}^{b}
&=\frac{1}{2}(\tau^*-\tau)g_{ij}.\\
\end{split}
\end{equation}
\noindent Similarly, setting $\lambda=\varDelta$ and $\mu=j$ in
\eqref{43}, and taking account of \eqref{219} and \eqref{221}, we
have
\begin{equation}\label{43b}
\begin{split}
\frac{1}{2}\xi^a{R_{jac}}^b{\phi_b}^c=\xi^a{\phi_j}^b\rho_{ab}.
\end{split}
\end{equation}
\noindent Transvecting \eqref{43a} with $g^{ij}$, and taking account
of (6) of \eqref{35}, we have
\begin{equation*}
\begin{split}
2\tau^*-2\tau+2-trh^2=\frac{3}{2}(\tau^*-\tau),
\end{split}
\end{equation*}
and hence
\begin{equation}\label{43c}
\begin{split}
\tau^*-\tau+4=2trh^2.
\end{split}
\end{equation}
\noindent Since a contact metric manifold is K-contact if and only
if $h=0$ holds (\cite{B2}, Theorem 7.1), from \eqref{43c}, taking
account of the Olszak's Lemma introduced in Section 3, we may obtain
the well-known result that a 3-dimensional K-contact manifold is
Sasakian. Then, taking account of Proposition \ref{Pr2}, we can
easily check that the equality \eqref{43b} is automatically
satisfied, since $\rho_{ia}\xi^a = 2\eta_i$ holds.\\
On the other hand, it is well-known that the curvature tensor $R$ of
a 3-dimensional Riemannian manifold $M=(M, g)$ satisfies the
following identity:
\begin{equation}\label{43d}
\begin{split}
R(X, Y)Z=&g(Y, Z)QX-g(X, Z)QY+g(QY, Z)X-g(QX, Z)Y\\
&-\frac{\tau}{2}(g(Y, Z)X-g(X, Z)Y),
\end{split}
\end{equation}
for any ${X}, {Y}, {Z} \in\mathfrak{X}({M})$. %Thus, from \eqref{24}
We here assume that $M=(M, \phi, \xi, \eta, g)$ is a $H$-contact
manifold. Then, there exists a smooth function $f$ on $M$ satisfying
\begin{equation}\label{45b}
Q\xi = f\xi
\end{equation}
\cite{P}. From \eqref{45b}, taking account of (5) of \eqref{35}, we
have
\begin{equation}\label{410}
\begin{split}
f=2-trh^2.
\end{split}
\end{equation}
Now from \eqref{43d} and \eqref{45b}, we have
\begin{equation}\label{411}
\begin{split}
{R_{ijk}}^l\xi^k=(\frac{\tau}{2}-f)(\eta_i{\delta_j}^l-\eta_j{\delta_i}^l)-\eta_i{\rho_j}^l+\eta_j{\rho_i}^l.
\end{split}
\end{equation}
From \eqref{411}, we have
\begin{equation*}
\begin{split}
{R_{ijk}}^l\xi^j\xi^k=(\frac{\tau}{2}-2f)\eta_i\xi^l-(\frac{\tau}{2}-f){\delta_i}^l+{\rho_i}^l,
\end{split}
\end{equation*}
and hence
\begin{equation}\label{412}
\begin{split}
{\phi_a}^i{R_{ijk}}^a\xi^j\xi^k=-(\frac{\tau}{2}-f){\phi_a}^l+{\phi_a}^i{\rho_i}^l.
\end{split}
\end{equation}
Operating $\nabla^i$ to \eqref{411} and taking account of (3) of
\eqref{35}, we have
\begin{equation}\label{413}
\begin{split}
&-\xi^k\nabla_k{\rho_j}^l+\xi^k\nabla^l\rho_{jk}+{R_{jik}}^l\phi^{ik}+{R_{jik}}^l{\phi_a}^kh^{ia}\\
&=(\frac{1}{2}\xi\tau-\xi{f}){\delta_j}^l+(\frac{\tau}{2}-f)(-{\phi_j}^l+\phi_{aj}h^{la})-{\rho_i}^l(-{\phi_j}^i+\phi_{aj}h^{ia}).
\end{split}
\end{equation}
Thus, from \eqref{413}, we have further
\begin{equation*}
\begin{split}
-\xi\tau+\frac{1}{2}\xi\tau-{h_a}^k{\phi_i}^a{\rho_k}^i=\frac{3}{2}\xi\tau-3\xi{f}+{h_a}^i{\phi_l}^a{\rho_i}^l,
\end{split}
\end{equation*}
and hence
\begin{equation}\label{414}
\begin{split}
2\xi\tau-3\xi{f}+2tr(h\phi{Q})=0.
\end{split}
\end{equation}
On one hand, from (4) of \eqref{35} and \eqref{412}, we have also
\begin{equation}\label{415}
\begin{split}
\xi^i\nabla_i{h_k}^l={\phi_k}^l-{h_a}^l{h_b}^a{\phi_k}^b+(\frac{\tau}{2}-f){\phi_k}^l-{\phi_a}^l{\rho_k}^a.
\end{split}
\end{equation}
Thus from \eqref{415}, taking account of \eqref{29}, we have
\begin{equation*}
\begin{split}
\frac{1}{2}\xi(trh^2)={h_l}^k\xi^i\nabla_i{h_k}^l=-{h_l}^k{h_a}^l{h_b}^a{\phi_k}^b-{h_l}^k{\phi_a}^l{\rho_k}^a=-tr(h^3\phi)-tr(h\phi{Q})=-tr(h\phi{Q}),
\end{split}
\end{equation*}
and hence
\begin{equation}\label{416}
\begin{split}
\frac{1}{2}\xi(trh^2)=-tr(h\phi{Q}).
\end{split}
\end{equation}
Thus, from \eqref{414} and \eqref{416}, taking account of
\eqref{410}, we have
\begin{equation*}
\begin{split}
0=2\xi\tau-3\xi{f}-\xi(trh^2)=2\xi\tau+3\xi(trh^2)-\xi(trh^2)=2\xi\tau+2\xi(trh^2),
\end{split}
\end{equation*}
and hence
\begin{equation}\label{417}
\begin{split}
\xi(\tau+trh^2)=0.
\end{split}
\end{equation}
We here recall the fact that there exists a 3-dimensional $(\kappa,
\mu, \nu)$-contact metric manifold $M_0$ which is open and close
\cite{K3}. We may easily check that the equality \eqref{417} is
regarded as an extension of the one on $M_0$ to the whole of $M$
which is derived from \eqref{320} and (3.21) with \eqref{313} in the
paper \cite{JPS} and also that the examples introduced in the same
paper illustrate the equality \eqref{417}.\\

\noindent {\bf Acknowledgements.} This work was supported by the
National Research Foundation of Korea (NRF) grant funded by the
Korea government (MEST) (2013020825).
%%%%%%%%%%%%%%%%%%%%%%%%%%%%%%%%%%%%%%%%%%%%%%%%%%%%%%%%%%%%%%%%%%%%%%%%%%%%%%%%%%%%%%%

\medskip

\noindent

\end{document}